# A CHARACTERIZATION OF HEDGING PORTFOLIOS FOR INTEREST RATE CONTINGENT CLAIMS

By Rene Carmona and Michael Tehranchi

*Princeton University and University of Texas at Austin*

We consider the problem of hedging a European interest rate contingent claim with a portfolio of zero-coupon bonds and show that an HJM type Markovian model driven by an infinite number of sources of randomness does not have some of the shortcomings found in the classical finite-factor models. Indeed, under natural conditions on the model, we find that there exists a unique hedging strategy, and that this strategy has the desirable property that at all times it consists of bonds with maturities that are less than or equal to the longest maturity of the bonds underlying the claim.

**1. Introduction.** This paper seeks to characterize portfolios that hedge contingent claims in the fixed income market. The fundamental traded instruments in this market are (zero-coupon) bonds, contracts in which the issuer agrees to pay one unit of currency at a fixed future maturity date. The idealized bonds considered here do not suffer from credit risk, that is, at maturity the bond issuer always makes the promised payment.

There are bonds with so many maturity dates traded on the market, it is conventional to assume at every time $t \geq 0$ there exists a bond that matures at time $T$ for every $T \geq t$. We use the notation $P_t(T)$ to denote the price at time $t$ of a bond with maturity date $T$.

Assuming that there is a continuum of traded securities is an important distinction from the classical Black–Scholes theory. Indeed, whereas in the Black–Scholes setting we work with a finite-dimensional vector $(S_t^1, \ldots, S_t^n)$ of stock prices at time $t$, in the fixed income market we work with the infinite-dimensional vector of the bond price curve $P_t(\cdot)$. It comes as no surprise then that the characterization of hedging porfolios in the fixed income market is a more subtle problem.

In Section 2 we review the classical finite-factor HJM models for the dynamics of the term structure of interest rates and discuss one of their major shortcomings: they allow for unnatural hedging strategies which would never be used by traders. In particular, for HJM models driven by a $d$-dimensional Wiener process, every interest rate contingent claim can be hedged perfectly by a portfolio of bonds of $d$ arbitrary maturities chosen a priori independently of the contingent claim. This result is at odds with traders' intuition that the maturities of the hedging bonds should depend on the contingent claim in question, and it is the main motivation for the present work.

In Section 3 we consider the natural generalization of HJM models driven by an infinite-dimensional Wiener process. We introduce the necessary functional analysis notation, and we define the function spaces which we use as state spaces for our infinite-dimensional dynamics. In this setting, the appropriate notion of portfolio is not obvious. But if we agree to consider a certain class of portfolios, we show that if a contingent claim can be hedged by a given strategy, then this strategy is unique under an appropriate model assumption. However, we run into two technical difficulties. First, for a given contingent claim it is not clear whether a hedging strategy exists at all. Second, if a strategy does exist, it is not obvious if the strategy agrees with the traders' intuition.

We are able to resolve these two technical problems with the tools of Malliavin calculus. Section 4 reviews briefly some useful results from this theory, including an infinite-dimensional version of the original Clark–Ocone formula. The results of Section 4 are *essentially* known. We state them clearly, and we prove those we could not find in the existing literature in an appropriate form.

In Section 5 we present the main results of this article. We consider the problem of hedging a European contingent claim for an infinite-factor Markovian HJM model where the payout functional is assumed to be Lipschitz. We explicitly compute the hedging strategy via the Clark–Ocone formula and show that the difficulties of the finite-factor HJM models can be overcome. In particular, under natural conditions on the model, we find in Theorem 5.7 that there exists a unique hedging strategy with the intuitively appealing property that at all times it consists of bonds with maturities that are less than or equal to the longest maturity of the bonds underlying the claim.

**2. Shortcomings of the finite-factor HJM models.** An important class of models of the fixed income market, introduced by Heath, Jarrow and Morton (1992) and henceforth called HJM models, takes the forward rate curve as the fundamental object to model. We can define the forward rate $f_t(T)$ at time $t$ for maturity $T$ to be given by the formula

$$(1) \qquad f_t(T) = -\frac{\partial}{\partial T} \log P_t(T)$$

whenever the bond price function is differentiable. Since a dollar today is worth more than a dollar tomorrow, we note that the bond price function $P_t(\cdot)$ is decreasing, implying by Vitali's theorem that the forward rate $f_t(T)$ exists for Lebesgue almost every $T \in [t, \infty)$. We assume in fact that the forward rates exist for every $T$, and in particular we can define the short rate $r_t$ at time $t$ by the relation

$$(2) \qquad r_t = f_t(t).$$

Note that the functions $P_t(\cdot)$ and $f_t(\cdot)$ contain the same information because the bond prices can be recovered from the forward rates via the equation

$$(3) \qquad P_t(T) = \exp\left(-\int_t^T f_t(s)\,ds\right).$$

The classical HJM model is specified by fixing a measurable space $(\Omega, \mathcal{F})$ and a risk-neutral probability measure $\mathbb{Q}$ for which there exists a standard $d$-dimensional Wiener process $\{W_t = (W_t^1, \ldots, W_t^d)\}_{t \geq 0}$ and the filtration $\{\mathcal{F}_t\}_{t \geq 0}$ given by the augmentation of the filtration generated by the Wiener process, such that the dynamics of the forward rate processes $\{f_t(T)\}_{t \in [0,T]}$ are given by

$$(4) \qquad df_t(T) = \left\langle \tau_t(T), \int_t^T \tau_t(s)\,ds \right\rangle_{\mathbb{R}^d} dt - \langle \tau_t(T), dW_t \rangle_{\mathbb{R}^d}$$

where $\{\tau_t(T) = (\tau_t^1(T), \ldots, \tau_t^d(T))\}_{t \in [0,T]}$ is an $\mathbb{R}^d$-valued adapted process for each $T > 0$ and the bracket $\langle \cdot, \cdot \rangle_{\mathbb{R}^d}$ is the usual Euclidean scalar product. The specific form of the drift was shown in Heath, Jarrow and Morton (1992) to be necessary to prohibit arbitrage. An important feature of this methodology is that the initial condition for these models is the whole forward rate curve $f_0(\cdot)$.

Note that for these models there are an infinite number of stochastic differential equations, one for each value of $T$, driven by a finite number of sources of randomness. Besides the fact that finite-dimensional Wiener processes are mathematically easier to handle than infinite-dimensional ones, this modeling assumption is usually justified by appealing to the statistics of the yield and forward rate curves observed on the market. [Recall that the yield $y_t(T)$ at time $t$ for maturity $T$ is given by $y_t(T) = (T-t)^{-1} \int_t^T f_t(s)\,ds$.] The principal component analysis of the U.S. Treasury yield curve as reported by Litterman and Scheinkman (1991) and of the Eurodollar forward rates by Bouchaud, Cont, Karoni, Potters and Sagna (1999) suggests that the dynamics of the forward rate are driven by a few sources of noise. Indeed,

Litterman and Scheinkman found that over 95% of the variations of the yield curve can be attributed to the first three factors, lending credence to HJM models driven by a low-dimensional Wiener process.

However, the assumption that the driving noise is finite dimensional has an important implication. Consider the problem of replicating the real $\mathcal{F}_T$-measurable random variable $\xi$ corresponding to the payout of an interest rate contingent claim that matures at a time $T$. Our hedging instruments are naturally the set of zero-coupon bonds and the risk-free bank account process $\{B_t\}_{t\geq 0}$ defined by

$$B_t = \exp\left(\int_0^t r_s\,ds\right). \tag{5}$$

To ease notation, we begin with a definition.

DEFINITION 2.1. For every process $\{X_t\}_{t\geq 0}$, we define the discounted process $\{\tilde{X}_t\}_{t\geq 0}$ by $\tilde{X}_t = B_t^{-1} X_t$. For every $\mathcal{F}_T$-measurable random variable $\xi$ corresponding to the payout of a contingent claim with maturity $T$, we use the notation $\tilde{\xi} = B_T^{-1}\xi$.

A major shortcoming of the finite-factor models is found in the following well-known proposition. We give a complete proof of this result to emphasize the difficulties we have to overcome in order to resolve the issues it raises.

PROPOSITION 2.2. *Suppose there exist $d$ dates $T_1 < T_2 < \cdots < T_d$ and a positive constant $c$ such that, for all $T < T_1$, the $d \times d$ matrix*

$$\sigma_t = \left[\tilde{P}_t(T_i)\int_t^{T_i} \tau_t^j(s)\,ds\right]_{i,j=1,\ldots,d} \tag{6}$$

*satisfies $\|\sigma_t x\|_{\mathbb{R}^d} \geq c\|x\|_{\mathbb{R}^d}$ for all $x \in \mathbb{R}^d$ and almost all $(t,\omega) \in [0,T] \times \Omega$. Then for every contingent claim $\xi$ with maturity $T < T_1$ such that $\mathbb{E}\{\tilde{\xi}^2\} < +\infty$, there exists a replicating strategy consisting of bonds with maturities $T_1, T_2, \ldots, T_d$ and the bank account.*

PROOF. Consider a strategy such that, at time $t$, the portfolio consists of $\phi_t^i$ units of the bond with maturity $T_i$ for $i=1,\ldots,d$ and of $\psi_t$ units of the bank account. As usual, we insist that our wealth process $\{V_t = \langle \phi_t, P_t\rangle_{\mathbb{R}^d} + \psi_t B_t\}_{t\geq 0}$ satisfies the self-financing condition

$$dV_t = \langle \phi_t, dP_t\rangle_{\mathbb{R}^d} + \psi_t\,dB_t,$$

where $\phi_t = (\phi_t^1,\ldots,\phi_t^d)$ is the vector of portfolio weights and $P_t = (P_t(T_1), \ldots, P_t(T_d))$ is the vector of bond prices. We now show that there exist processes $\{\phi_t\}_{t\in[0,T]}$ and $\{\psi_t\}_{t\in[0,T]}$ such that $V_T = \xi$ almost surely.

Recall that the bond price at time $t$ for maturity $s$ is related to the forward rates by (3), so by an application of Itô's rule and the stochastic Fubini's theorem we have that the dynamics of the bond price for each $T_i$ satisfy the equation

$$\frac{dP_t(T_i)}{P_t(T_i)} = r_t\, dt + \left\langle \int_t^{T_i} \tau_t(u)\, du,\, dW_t \right\rangle_{\mathbb{R}^d}. \tag{7}$$

By (7) and (6), the dynamics of the vector of discounted bond prices is given by $d\tilde{P}_t = \sigma_t\, dW_t$, and consequently, the dynamics of the discounted wealth process is given by

$$d\tilde{V}_t = \langle \phi_t, d\tilde{P}_t \rangle_{\mathbb{R}^d} = \langle \sigma_t^* \phi_t, dW_t \rangle_{\mathbb{R}^d}.$$

On the other hand, if $\mathbb{E}\{\tilde{\xi}^2\} < +\infty$, we can apply Itô's martingale representation theorem to conclude that there exists a $d$-dimensional adapted process $\{\alpha_t\}_{t\in[0,T]}$ such that $\mathbb{E}\{\int_0^T \|\alpha_t\|_{\mathbb{R}^d}^2\, dt\} < +\infty$ and

$$\tilde{\xi} = \mathbb{E}\{\tilde{\xi}\} + \int_0^T \langle \alpha_t, dW_t \rangle_{\mathbb{R}^d}.$$

Setting the initial wealth $V_0 = \mathbb{E}\{\tilde{\xi}\}$ and portfolio weights $\phi_t = \sigma_t^{*-1} \alpha_t$ and $\psi_t = \tilde{V}_t - \langle \phi_t, \tilde{P}_t \rangle_{\mathbb{R}^d}$, we find our desired replicating strategy. $\square$

Thus, in this $d$-factor HJM model every square-integrable claim can be replicated by a strategy of holding bonds maturing at the $d$ dates $T_1, \ldots, T_d$ fixed a priori and independently of the claim. For instance, with a three-factor HJM model, it is possible to perfectly hedge a call option on a bond of maturity five years with a portfolio of bonds of maturity fifteen, twenty and twenty-five years. Cont (2004) remarks that this result is counterintuitive and contrary to market practice. Indeed, there seems to be a notion of "maturity specific risk" not captured by finite-factor HJM models since we expect that such a contingent claim should be hedged with bonds of maturities less than or equal to five years. This shortcoming can be attributed to the high degree of redundancy in the finite-factor models.

In Section 3 we show that if the dynamics of the bond prices are driven by an infinite-dimensional Wiener process, we can find conditions on the model such that a given hedging strategy is unique. Unfortunately, the usual notions of hedging become more complicated in infinite dimensions.

**3. Infinite-factor HJM models: some difficulties.** In this section we take a first look at the hedging problem for infinite-factor HJM models, our goal being to emphasize some of the difficulties occurring because of the infinite-dimensionality of the sources of randomness.

Here and throughout the rest of the paper, the stochastic processes are assumed to be defined on a complete probability space $(\Omega, \mathcal{F}, \mathbb{Q})$. Also, for the ease of exposition we prefer to break from the HJM tradition and choose the state variable for these models to be the discounted bond price curve $\tilde{P}_t(\cdot)$ instead of the forward rate curve $f_t(\cdot)$. But noting that the price of a bond at maturity is $P_t(t) = 1$, we see that the bank account process can be recovered by the formula

$$B_t = \frac{1}{\tilde{P}_t(t)} \tag{8}$$

and the bond price with maturity $s \geq t$ can be recovered via

$$P_t(s) = \frac{\tilde{P}_t(s)}{\tilde{P}_t(t)}. \tag{9}$$

This change of variables eases the analysis, although it is quite superficial in the sense that there is a one-to-one correspondence between bond prices and instantaneous forward rates given by (1). As motivation for this change of variables, consider a European option that matures at time $T$ and has a payout of the form

$$\xi = g(P_T(T_1), \ldots, P_T(T_n)) \tag{10}$$

for some dates $T_i \geq T$ and some measurable function $g : \mathbb{R}^n \to \mathbb{R}$. From the discussion in Section 2 we see that, in order to replicate such a claim, we must find the martingale representation of the discounted claim $\tilde{\xi} = B_T^{-1} g(P_T(T_1), \ldots, P_T(T_n))$. By (5) we see that $\tilde{\xi}$ depends not only on the bond prices at time $T$ but also on the entire history of the short rate process. But treating the discounted bond price curve as the state variable, we have by (8) and (9) that

$$\tilde{\xi} = \tilde{P}_T(T) \, g\left(\frac{\tilde{P}_T(T_1)}{\tilde{P}_T(T)}, \ldots, \frac{\tilde{P}_T(T_n)}{\tilde{P}_T(T)}\right).$$

Defining the functional $\tilde{g}$ on the space $C(\mathbb{R}_+)$ of continuous functions on $\mathbb{R}_+$ by

$$\tilde{g}(x) = x(T) g\left(\frac{x(T_1)}{x(T)}, \ldots, \frac{x(T_n)}{x(T)}\right), \tag{11}$$

we have that $\tilde{\xi} = \tilde{g}(\tilde{P}_T)$ only depends on the time $T$ values of the discounted bond price processes. Of course this trades one infinite-dimensional problem for another; yet in this framework, the problem can be treated as the Black–Scholes problem of pricing and hedging a modified European contingent claim on a portfolio of "stocks" with zero interest rate.

REMARK 3.1 (Settlement in arrears). There are many interest rate options that pay in arrears. That is, although the payout $\xi = g(P_T(T_1), \ldots, P_T(T_n))$ is $\mathcal{F}_T$-measurable, the money does not change hands until the future settlement date $T + \Delta T$. This is the case for claims contingent on the LIBOR rate, such as caplets. In this situation the discounted claim is

$$\tilde{\xi} = B_{T+\Delta T}^{-1} g(P_T(T_1), \ldots, P_T(T_n)).$$

Noting that $\tilde{P}_{T+\Delta T}(T+\Delta T) = \tilde{P}_T(T+\Delta T) + \int_T^{T+\Delta T} d\tilde{P}_t(T+\Delta T)$, we have by (8) and (9) that

$$\tilde{\xi} = \hat{g}(\tilde{P}_T) + \int_T^{T+\Delta T} \xi \, d\tilde{P}_t(T+\Delta T),$$

where $\hat{g} : C(\mathbb{R}_+) \to \mathbb{R}$ is defined by

$$\hat{g}(x) = x(T+\Delta T) g\left(\frac{x(T_1)}{x(T)}, \ldots, \frac{x(T_n)}{x(T)}\right).$$

Thus the strategy that consists of replicating the $\mathcal{F}_T$-measurable random variable $\hat{g}(\tilde{P}_T)$ and then holding $\xi$ units of the bond with maturity $T + \Delta T$ replicates the payout of the contingent claim. Hence the hedging problem still maps to an infinite-dimensional zero interest rate Black–Scholes world, but with the payout function modified slightly differently.

REMARK 3.2. For each time $t$, the domain of the discounted bond price function $\tilde{P}_t(\cdot)$ is the interval $[t, \infty)$. Since we need to consider the dynamics of the discounted price curve as $t$ varies, it would be more convenient if the functions $\tilde{P}_t(\cdot)$ had a common domain. For this reason, we assume that, for every $t \geq 0$, the domain of the discounted bond price function $\tilde{P}_t(\cdot)$ is the interval $[0, \infty)$, where we extend the definition of $\tilde{P}_t(\cdot)$ by $\tilde{P}_t(s) = B_s^{-1}$ for $s \in [0, t]$. Note then that the process $\{\tilde{P}_t(s)\}_{t \geq 0}$ is constant for $t \geq s$. The corresponding bond prices are given by $P_t(s) = B_s^{-1} B_t$ for $s \in [0, t]$ so that this extension conforms with the price

(12) $$P_t(s) = \mathbb{E}\left\{\frac{B_t}{B_s} \,\Big|\, \mathcal{F}_t\right\}$$

for $s \geq t$ and can be understood that once a bond matures the one dollar payout is immediately put into the bank to accrue interest at the short rate.

REMARK 3.3. Another popular way to resolve the issue of having functions defined on time-dependent domains is to switch to the so-called Musiela notation. The idea is to work with the time *to* maturity $\theta = T - t$ rather than with the time *of* maturity $T$. In this approach, the reparameterized discounted bond price curve $\hat{P}_t(\cdot)$ is defined by

$$\hat{P}_t(\theta) = \tilde{P}_t(t + \theta).$$

For the finite-factor HJM model, this new process is a weak solution of the following stochastic partial differential equation:

$$d\hat{P}_t(\theta) = \frac{\partial \hat{P}_t(\theta)}{\partial \theta}\,dt + \hat{P}_t(\theta)\Big\langle \int_t^{t+\theta} \tau_t(s)\,ds, dW_t \Big\rangle_{\mathbb{R}^d}.$$

This formulation of the HJM models proved to be very fruitful; see, for instance, Musiela (1993), Goldys and Musiela (1996) and Filipovic (2001). Nevertheless, for the sake of studying the hedging strategies for interest rate contingent claims, it is more convenient to retain the time to maturity parameterization. Indeed, whereas the process $\{\tilde{P}_t(T)\}_{t\geq 0}$ is a martingale for each $T$, the analogous process in Musiela notation $\{\hat{P}_t(\theta)_t\}_{t\geq 0}$ is usually not a martingale for any $\theta$.

In order to specify an infinite-factor model of the evolution of the discounted bond prices, it is natural to work in a function space setting. We first review the relevant notation of functional analysis. For a Banach space $E$, the duality form is denoted $\langle \cdot, \cdot \rangle_E : E^* \times E \to \mathbb{R}$. If $F$ is another Banach space, we let $\mathcal{L}(F, E)$ denote the Banach space of bounded linear operators taking $F$ into $E$ with norm

$$\|A\|_{\mathcal{L}(F,E)} = \sup_{x \in F,\, \|x\|_F \leq 1} \|Ax\|_E.$$

If $A \in \mathcal{L}(F, E)$, the (Banach space) adjoint $A^*$ of $A$ is the unique element of $\mathcal{L}(E^*, F^*)$ satisfying

$$\langle \mu, Ax \rangle_E = \langle A^*\mu, x \rangle_F \qquad \text{for all } \mu \in E^*,\ x \in F.$$

If $G$ is a Hilbert space, we use the notation $x^* \in G^*$ for the Riesz representation of the element $x \in G$, and we identify the double dual $G^{**}$ with $G$. If $\mathcal{S}$ is a subspace of $G$, then we let

$$\mathcal{S}^\perp = \{\mu \in G^* \text{ such that } \langle \mu, x \rangle_G = 0 \text{ for all } x \in G\}$$

be the closed subspace of $G^*$ orthogonal to $\mathcal{S}$.

If $G$ and $H$ are separable Hilbert spaces, the space of Hilbert–Schmidt operators taking $H$ into $G$ is denoted $\mathcal{L}_{\mathrm{HS}}(H, G)$ and is itself a Hilbert space for the norm

$$\|A\|_{\mathcal{L}_{\mathrm{HS}}(H,G)} = \left( \sum_{i=1}^{\infty} \|Ae_i\|_G^2 \right)^{1/2},$$

where $\{e_i\}_i$ is any orthonormal basis for $H$. There is a natural isometry of the space $\mathcal{L}_{\mathrm{HS}}(H, G)$ and the Hilbert space tensor product $G \otimes H^*$.

For a Banach space $E$, we denote by $L^p(\Gamma; E)$ the Banach space of (equivalence classes of) measurable functions from $\Gamma$ into $E$ with the norm

$$\|f\|_{L^p(\Gamma;E)} = \left(\int_\Gamma \|f(x)\|_E^p \mu(dx)\right)^{1/p},$$

where the measure space $(\Gamma, \mathcal{G}, \mu)$ is the interval $([0,T], \mathcal{B}_{[0,T]}, \text{Leb}_{[0,T]})$, the probability space $(\Omega, \mathcal{F}, \mathbb{Q})$, or their product $([0,T] \times \Omega, \mathcal{B}_{[0,T]} \otimes \mathcal{F}, \text{Leb}_{[0,T]} \times \mathbb{Q})$.

For our application, we need an infinite-dimensional version of the vector-valued stochastic integrals of the form $\int_0^t \sigma_s \, dW_s$. Self-contained expositions of the theory of infinite-dimensional stochastic integration can be found in the books of Da Prato and Zabczyk (1992), Kallianpur and Xiong (1995) and Carmona (2004). From now on, we fix a real separable Hilbert space $H$, and we assume that $\{W_t\}_{t \geq 0}$ is a cylindrical $H$-valued Wiener process defined on the probability space $(\Omega, \mathcal{F}, \mathbb{Q})$, that this cylindrical process generates the $\sigma$-field $\mathcal{F}$, and that the filtration $\{\mathcal{F}_t\}_{t \geq 0}$ is given by the augmentation of the filtration it generates. The classical finite-factor HJM model corresponds to the choice of a finite-dimensional space $H$. The integrands considered here are the adapted, measurable and square-integrable stochastic processes $\sigma = \{\sigma_t\}_{t \geq 0}$ valued in the space $\mathcal{L}_{\text{HS}}(H, F)$ of Hilbert–Schmidt operators from $H$ into $F$ for which we have Itô's isometry

$$\mathbb{E}\left\{\left\|\int_0^t \sigma_s \, dW_s\right\|_F^2\right\} = \mathbb{E}\left\{\int_0^t \|\sigma_s\|_{\mathcal{L}_{\text{HS}}(H,F)}^2 \, ds\right\}.$$

Note if $F = \mathbb{R}$, the space $\mathcal{L}_{\text{HS}}(H, \mathbb{R})$ of Hilbert–Schmidt operators is just $H^*$. In this case we write

$$\int_0^t \sigma_s \, dW_s = \int_0^t \langle \sigma_s, dW_s \rangle_H$$

in analogy with the finite-dimensional stochastic integration. However, this notation can only be formal if $H$ is infinite dimensional since the Wiener process $\{W_t\}_{t \geq 0}$ visits the space $H$ with probability zero.

We now introduce a family of weighted Sobolev spaces to serve as the state space for the infinite-dimensional dynamics.

DEFINITION 3.1. For every function $w : \mathbb{R}_+ \to \mathbb{R}_+$ and for $i = 1, 2$, we define the space $F_w^i$ of functions $x : \mathbb{R}_+ \to \mathbb{R}$ which are $i-1$ times differentiable, with a $(i-1)$st derivative $x^{(i-1)}$ absolutely continuous and such that $x^{(j)}(\infty) = 0$ for $j \leq i-1$, and $\int_0^\infty x^{(i)}(u)^2 w(u) \, du < +\infty$.

The space $F_w^i$ is a Hilbert space for the norm $\|x\|_{F_w^i} = (\int_0^\infty x^{(i)}(u)^2 w(u) \, du)^{1/2}$.

We work with the spaces $F_v^1$ and $F_w^2$. We list three useful properties of these spaces.

PROPOSITION 3.2.   *If the positive function $v$ is such that $C_v = \int_0^\infty v(s)^{-1} \, ds < +\infty$, then the evaluation functionals $\delta_s$, where $\langle \delta_s, x \rangle_{F_v^1} = x(s)$, are continuous on $F_v^1$ for all $s \geq 0$.*

*If the positive function $w$ is such that $C_w = \int_0^\infty (1 + u^2) w(u)^{-1} \, du < +\infty$, then the evaluation functionals $\delta_s$ and the point-wise differentiation $\delta_s'$, where $\langle \delta_s', x \rangle_{F_w^2} = -x'(s)$, are both linear continuous functionals on $F_w^2$ for all $s \geq 0$.*

*If $C_{vw} = \int_0^\infty \int_0^s v(u)/w(s) \, du \, ds < +\infty$, then the inclusion from $F_w^2$ to $F_v^1$ is continuous.*

PROOF.   The evaluation functionals are uniformly bounded on $F_v^1$ since

$$|\langle \delta_s, x \rangle_{F_v^1}| = |x(s)| = \left| \int_s^\infty x'(u) \, du \right|$$

$$\leq \left( \int_s^\infty \frac{du}{w(u)} \right)^{1/2} \left( \int_s^\infty x'(u)^2 w(u) \, du \right)^{1/2}$$

$$\leq C_v^{1/2} \|x\|_{F_v^1}.$$

Similarly, the point-wise differentiation functionals are uniformly bounded on $F_w^2$ by $C_v^{1/2}$. The evaluation functionals are uniformly bounded on $F_w^2$ since

$$|\langle \delta_s, x \rangle_{F_w^2}| = |x(s)| = \left| \int_s^\infty x'(u) \, du \right|$$

$$= \left| \int_s^\infty (u - s) x''(u) \, du \right|$$

$$\leq \left( \int_s^\infty \frac{(u-s)^2 \, du}{w(u)} \right)^{1/2} \left( \int_s^\infty x''(u)^2 w(u) \, du \right)^{1/2}$$

$$\leq C_w^{1/2} \|x\|_{F_w^2}.$$

Finally, if $x \in F_v^1 \cap F_w^2$, then

$$\|x\|_{F_v^1}^2 = \int_0^\infty x'(s)^2 v(s) \, ds$$

$$\leq \int_0^\infty \int_s^\infty \frac{du}{w(u)} \|x\|_{F_w^1}^2 v(s) \, ds = C_{vw} \|x\|_{F_w^1}^2,$$

and hence the inclusion from $F_w^2$ to $F_v^1$ is continuous.   □

We fix a weight $w$ satisfying the conditions of Proposition 3.2, and from now on we assume that the state variable $\tilde{P}_t(\cdot)$ is an element of the function space $F_w^2$ for every $t \geq 0$. Note that for this choice of state space, we may speak honestly about the price of a specific bond or the value of a specific forward rate since evaluation and point-wise differentiation are continuous. This choice also agrees with the fact that a bond that never matures is worthless and hence $\tilde{P}_t(\infty) = B_t^{-1} P_t(\infty) = 0$. And since $F_w^2$ is a Hilbert space, we may use the integration theory mentioned above.

We now formulate a model of the discounted price dynamics.

ASSUMPTION 3.3. The risk-neutral dynamics of the discounted price curve $\{\tilde{P}_t\}_{t \geq 0}$ are described by the initial condition $\tilde{P}_0 \in F_w^2$ and the evolution equation

$$(13) \qquad d\tilde{P}_t = \sigma_t \, dW_t,$$

where $\{\sigma_t\}_{t \geq 0}$ is an $\mathcal{L}_{\mathrm{HS}}(H, F_w^2)$-valued adapted stochastic process such that

$$(14) \qquad \sigma_t^* \delta_s = 0 \qquad \text{for all } t \geq s.$$

We assume that $\tilde{P}_0$ and $\{\sigma_t\}_{t>0}$ conspire in such a way that $\tilde{P}_t(s) > 0$ for all $s \geq 0$ and that

$$(15) \quad \mathbb{E}\left\{ \int_0^t \left( \frac{|\tilde{P}_s'(s)| \, \|\tilde{P}_s\|_{F_w^2}}{\tilde{P}_s(s)^2} + \left(1 + \frac{1}{\tilde{P}_s(s)^2}\right) \|\sigma_s\|_{\mathcal{L}_{\mathrm{HS}}(H, F_w^2)}^2 \right) ds \right\} < +\infty.$$

REMARK 3.4. Condition (14) guarantees that the process $\{\tilde{P}_t(s)\}_{t \geq 0}$ becomes constant for $t \geq s$. Indeed, the continuity of $\delta_s$ implies that, for $t \geq s$, we have

$$\tilde{P}_t(s) - \tilde{P}_s(s) = \int_s^t \langle \sigma_u^* \delta_s, dW_u \rangle_H = 0.$$

Since our starting point is the discounted bound curve $\tilde{P}_t$, we need to infer the definition of the bank account $B_t$ and of the zero-coupon curve $P_t$. The bank account, given by the formula $B_t = \tilde{P}_t(t)^{-1}$, has dynamics given formally by

$$(16) \qquad dB_t = -\frac{\tilde{P}_t'(t)}{\tilde{P}_t(t)^2} \, dt,$$

while the prices of the zero-coupon bonds, given by $P_t = \tilde{P}_t / \tilde{P}_t(t)$, have dynamics given formally by

$$(17) \qquad dP_t = -\frac{\tilde{P}_t'(t)}{\tilde{P}_t(t)^2} \tilde{P}_t \, dt + \frac{1}{\tilde{P}_t(t)} \sigma_t \, dW_t.$$

Condition (15) ensures that the stochastic equations (13), (16) and (17) make perfectly good sense. In Remark 5.1 we give sufficient conditions for (15) to hold for a Markovian HJM model.

REMARK 3.5. Because the eigenvalues of a Hilbert–Schmidt operator must decay fast enough for the sum of their squares to be finite, assuming that $H$ is infinite dimensional does not disagree with the principal component analysis typically used to justify the introduction of models with finitely many factors.

Given such a model, we propose to study how to hedge a contingent claim. Equivalently, this problem is equivalent to the search for a representation of the contingent claim as a stochastic integral with respect to the underlying price process. As we are about to see, this task reduces to finding an adapted process $\{\phi_t\}_{t\geq 0}$ such that

$$\tilde{\xi} = \mathbb{E}\{\tilde{\xi}\} + \int_0^T \langle \phi_s, d\tilde{P}_s \rangle_{F_w^2}.$$

Identifying this process $\{\phi_t\}_{t\geq 0}$ illustrates the difficulties of working in infinite dimensions.

In the real world, a portfolio can only contain a finite number of bonds at any time. That is, we really should only consider processes such that, for almost every $(t,\omega) \in [0,T] \times \Omega$, we can find a positive integer $d$, positive real numbers $T_1, \ldots, T_d$ and real numbers $c_1, \ldots, c_d$ so that

$$\phi_t = \sum_{i=1}^d c_i \delta_{T_i}.$$

However, limiting ourselves to such portfolios at this stage of the analysis would be severely restrictive. Indeed, since we are willing to assume that there exists a continuum of traded securities, it seems reasonable to assume that we can form portfolios with bonds of an infinite number of different maturities. Since the process $\{\tilde{P}_t\}_{t\geq 0}$ takes values in $F_w^2$, it would seem natural to require that $\{\phi_t\}_{t\in[0,T]}$ takes values in the dual $F_w^{2*}$. Remember that measures of the form $\sum_{i=1}^d c_i \delta_{T_i}$ are in $F_w^{2*}$. But the elements of the space $F_w^2$ are functions that are quite smooth, and consequently, the dual space $F_w^{2*}$ contains distributions that can be quite rough. Indeed, point-wise differentiation is bounded on $F_w^2$, and we choose to work in this space precisely because we need to define the short rate in the drift of the bond price process. Nevertheless, even though we would prefer to think of our hedging strategies as being measures, if we work with $F_w^{2*}$-valued portfolios, we risk the uncomfortable possibility that they might be much wilder distributions.

As a partial resolution to this problem, we consider strategies valued in $F_v^{1*}$, where $v$ is a function satisfying the conditions of Proposition 3.2. Since the inclusion map from $F_w^2$ into $F_v^1$ is continuous, the dual $F_v^{1*}$ can be identified with a dense subset of $F_w^{2*}$. We fix a $v$ once and for all, and henceforth adopt the notation $F = F_v^1$.

We now make precise the various notions of strategy we shall use.

DEFINITION 3.4. A *strategy* is an adapted $F^*$-valued process $\{\phi_t\}_{t\geq 0}$ such that $\phi_t \in \mathcal{T}_{[t,\infty)}$ for almost every $(t,\omega)$, where we use the notation

$$\mathcal{T}_A = \overline{\operatorname{span}\{\delta_s; s \in A\}} \subset F^* \tag{18}$$

for a closed interval $A \subset \mathbb{R}_+$ and where the closure is taken in the topology of $F^*$.

Note that the restriction $\phi_t \in \mathcal{T}_{[t,\infty)}$ reflects the fact that it is unnecessary to hold expired bonds.

DEFINITION 3.5. A *self-financing strategy* is a strategy $\{\varphi_t\}_{t\geq 0}$ such that $d\langle \varphi_t, P_t \rangle_F = \langle \varphi_t, dP_t \rangle_F$.

For each strategy $\{\phi_t\}_{t\geq 0}$, the associated wealth process $\{V_t\}_{t\geq 0}$ has dynamics

$$dV_t = \langle \phi_t, dP_t \rangle_F + \psi_t \, dB_t$$

with $\psi_t B_t = V_t - \langle \phi_t, P_t \rangle_F$. But by (16) and (17), we have $dB_t = B_t \langle \delta_t, dP_t \rangle_F$; that is to say, the bank account can be replicated by the self-financing strategy of holding the bond maturing instantly. Hence, the dynamics of the wealth process can be written in the form

$$dV_t = \langle \phi_t + (V_t - \langle \phi_t, P_t \rangle_F)\delta_t, dP_t \rangle_F,$$

and for every strategy $\{\phi_t\}_{t\geq 0}$ we can construct a self-financing strategy $\{\varphi_t\}_{t\geq 0}$ via the rule

$$\varphi_t = \phi_t + (V_t - \langle \phi_t, P_t \rangle_F)\delta_t.$$

DEFINITION 3.6. A *pre-hedging strategy* for the contingent claim $\xi$ is a strategy $\{\phi_t\}_{t \in [0,T]}$ such that

$$\tilde{\xi} = \mathbb{E}\{\tilde{\xi}\} + \int_0^T \langle \phi_s, d\tilde{P}_s \rangle_F.$$

A *hedging strategy* for the contingent claim $\xi$ is a self-financing pre-hedging strategy.

Note that a pre-hedging strategy need not be self-financing since the pre-hedging condition is indifferent to the amount held in the bank account.

We now show that if a contingent claim can be hedged by an $F^*$-valued strategy, then under an appropriate model assumption, the hedging strategy is unique. This is a first step toward eliminating the counterintuitive strategies found in Section 2 in the case of finite-dimensional models.

Note that under the condition (14) we necessarily have $\ker(\sigma_t^*) \supset \mathcal{T}_{[0,t]}$ since $\sigma_t^*$ is almost surely a bounded operator. If we insist that this inclusion is an equality, we have the following proposition.

PROPOSITION 3.7. *Suppose, for almost all $(t, \omega) \in [0, T] \times \Omega$, we have*

(19) $$\ker(\sigma_t^*) = \mathcal{T}_{[0,t]}.$$

*If the hedging strategies $\{\varphi_t^1\}_{t \in [0,T]}$ and $\{\varphi_t^2\}_{t \in [0,T]}$ hedge the same claim $\xi$, then $\varphi_t^1 = \varphi_t^2$ for almost all $(t, \omega) \in [0, T] \times \Omega$.*

PROOF. Clearly the strategy $\{\varphi_t = \varphi_t^1 - \varphi_t^2\}_{t \in [0,T]}$ replicates the zero payout. Since $\{\varphi_t\}_{t \in [0,T]}$ is pre-hedging, we have $\int_0^T \langle \varphi_t, d\tilde{P}_t \rangle_F = 0$ almost surely and hence $\mathbb{E}\{\int_0^T \|\sigma_t^* \varphi_t\|_{H^*}^2 dt\} = 0$. Thus we have that $\varphi_t \in \mathcal{T}_{[0,t]} \cap \mathcal{T}_{[t,\infty)}$ for almost all $(t, \omega)$. Furthermore, since $\{\varphi_t\}_{t \in [0,T]}$ is self-financing, we have $\langle \varphi_t, P_t \rangle_F = 0$ for almost all $(t, \omega)$.

Fixing $(t, \omega)$, we now let $j \in F$ be any function with the property that $j(t) = 1$. Note that the function $\mathbf{1}_{[t,\infty)}(P_t - j)$ is in $F$ and that

$$\langle \varphi_t, \mathbf{1}_{[t,\infty)}(P_t - j) \rangle_F = 0.$$

Similarly, the function $\mathbf{1}_{[0,t]}(P_t - j)$ is also in $F$ and

$$\langle \varphi_t, \mathbf{1}_{[0,t]}(P_t - j) \rangle_F = 0.$$

Hence, we have $\langle \varphi_t, j \rangle_F = 0$ for every $j \in F$; thus $\varphi_t = 0$. □

REMARK 3.6. The above proof of uniqueness does not go through if we had allowed $F_w^{2*}$-valued portfolios. In particular, there exist nonzero portfolios $\varphi_t \in F_w^{2*}$ such that

$$\varphi_t \in \overline{\mathrm{span}\{\delta_s; s \leq t\}}^{F_w^{2*}} \cap \overline{\mathrm{span}\{\delta_s; s \geq t\}}^{F_w^{2*}}$$

and $\langle \varphi_t, P_t \rangle_F = 0$. For example, let $\varphi_t = r_t \delta_t - \delta_t'$. This is another reason for demanding that the portfolios be valued in the smaller space $F^* = F_v^{1*}$.

Only with infinite-dimensional $H$ can we hope to satisfy the conditions of the above proposition. However, unlike the finite-dimensional case, it is not clear that such hedging strategies exist in general. Parroting the calculation from Section 2, let $\tilde{\xi}$ be a square-integrable discounted claim and suppose that we could find a pre-hedging strategy $\{\phi_t\}_{t \geq 0}$ such that

$$\tilde{\xi} = \mathbb{E}\{\tilde{\xi}\} + \int_0^T \langle \phi_s, d\tilde{P}_s \rangle_F = \mathbb{E}\{\tilde{\xi}\} + \int_0^T \langle \sigma_s^* \phi_s, dW_s \rangle_H.$$

But recall that the martingale representation theorem states that there exists an adapted $H^*$-valued process $\{\alpha_t\}_{t \in [0,T]}$ such that $\mathbb{E} \int_0^T \|\alpha_s\|_{H^*}^2 ds < +\infty$ and

(20) $$\tilde{\xi} = \mathbb{E}\{\tilde{\xi}\} + \int_0^T \langle \alpha_s, dW_s \rangle_H.$$

See Da Prato and Zabczyk (1992) or Carmona (2004) for the infinite-dimensional version of this result. Thus, in order to calculate a pre-hedging portfolio at time $t$, we need only compute $\phi_t = \sigma_t^{*-1}\alpha_t$. But by assumption the operator $\sigma_t$ is Hilbert–Schmidt almost surely. Since $H$ is infinite dimensional, the inverse $\sigma_t^{*-1}$ is unbounded, and at this level of generality there is no guarantee that $\alpha_t$ is in its domain for any $t$. Thus restricting the portfolio to be in the space $F^*$ for all $t \geq 0$ is insufficient to replicate every square-integrable contingent claim. Björk, Di Masi, Kabanov and Runggaldier (1997) discuss this difficulty in the Banach space setting where the bond price process is a jump-diffusion driven by a finite-dimensional Wiener process, and they introduce the notion of approximate market completeness.

We could proceed by enlarging the class of allowable hedging portfolios by insisting that $\phi_t$ is in the so-called *covariance* space $\sigma_t^{*-1}H^*$ for almost all $(t,\omega) \in [0,T] \times \Omega$, where $\sigma^{*-1}H^* \supset F^*$ is the Hilbert space with norm $\|\phi\|_{\sigma^{*-1}H^*} = \|\sigma^*\phi\|_{H^*}$. De Donno and Pratelli (2004) elaborate on this approach for models in which the price process is defined cylindrically on a Hilbert space $F$. Notice the spaces $\sigma_t^{*-1}H^*$ generally depend on $t$ and $\omega$, but it would be nicer if the hedging strategy were valued in a fixed space with a more explicit characterization. Furthermore, we would need the bond price $P_t \in (\sigma_t^{*-1}H^*)^* = \sigma_t H$ to be in a much smaller space almost surely in order to construct the self-financing strategy.

Even if we knew that $\alpha_t$ was in the domain of $\sigma_t^{*-1}$, it would be unclear if the portfolio $\phi_t = \sigma_t^{*-1}\alpha_t$ agrees with the traders' intuition, since the support of $\phi_t$ is interpreted as the range of maturities of the bonds in the portfolio. We see that in order to construct a reasonable hedging portfolio, we need to know some detailed information about the martingale representation of the payout. In the classical Black–Scholes framework of a complete market with finitely many tradable assets, the hedging portfolio of a contingent claim is expressed as the gradient of the solution of a parabolic partial differential equation. Goldys and Musiela (1996) extend this PDE approach to the bond market setting by finding conditions under which the solution of the infinite-dimensional PDE is differentiable. In Section 5 we take a somewhat different approach to construct the hedging portfolio by appealing to the Clark–Ocone formula of Malliavin calculus.

Indeed, if we limit ourselves to payouts of the form $\xi = g(P_T)$, we can find conditions on the model paramaters $\{\sigma_t\}_{t \geq 0}$ and the payout function $g(\cdot)$ under which there exists a unique $F^*$-valued hedging strategy. Furthermore, under assumptions often satisfied by models used in practice, these conditions imply that the portfolio is confined to a small subspace of $F^*$.

**4. Malliavin calculus and the Clark–Ocone formula.** For an $\mathcal{F}_T$-measurable random variable $\tilde{\xi} \in L^2(\Omega; \mathbb{R})$, the martingale representation theorem guarantees the existence of an $H^*$-valued integrand such that $\tilde{\xi}$ can be written as a stochastic integral with respect to the Wiener process. For the financial application motivating this article, it is necessary to have an explicit formula for this integrand, expressed in terms of $\tilde{\xi}$. Fortunately, under a differentiability assumption on $\tilde{\xi}$, the Clark–Ocone formula provides such an expression. In order to state this useful result, we need to first introduce the Malliavin derivative operator and list some of its properties. The material of this section can be found in Nualart's (1995) book when the Hilbert spaces are finite dimensional.

The Malliavin derivative is a linear map from a space of random variables to a space of processes. We are concerned with the case where the random variables are elements of $L^2(\Omega; G)$, in which case the processes are elements $L^2([0, T] \times \Omega; \mathcal{L}_{\mathrm{HS}}(H, G))$, where $G$ is a real separable Hilbert space.

Being a derivative, it is not surprising that this operator is unbounded on $L^2(\Omega; G)$. We take the approach of defining it first on a core and then extending the definition to the closure of this set in the graph norm topology.

We now define the Malliavin derivative operator $D$ on this set.

DEFINITION 4.1. The random variables $X \in L^2(\Omega; G)$ of the form

$$(21) \qquad X = \kappa\bigg(\int_0^T \langle h_t^1, dW_t\rangle_H, \ldots, \int_0^T \langle h_t^n, dW_t\rangle_H\bigg),$$

where $h^1, \ldots, h^n \in L^2([0, T]; H^*)$ are deterministic, and where the differentiable function $\kappa : \mathbb{R}^n \to G$ is such that

$$(22) \qquad \sum_{i=1}^n \bigg\|\frac{\partial \kappa(x)}{\partial x_i}\bigg\|_G < C(1 + \|x\|_{\mathbb{R}^n}^p),$$

for some $p$, $C > 0$, and for all $x = (x_1, \ldots, x_n) \in \mathbb{R}^n$, are called *smooth*, and their Malliavin derivatives are defined to be

$$DX = \sum_{i=1}^n \frac{\partial \kappa}{\partial x_i}\bigg(\int_0^T \langle h_t^1, dW_t\rangle_H, \ldots, \int_0^T \langle h_t^n, dW_t\rangle_H\bigg) \otimes h^i.$$

Note that the process $\{D_t X\}_{t \in [0, T]}$ is valued in $\mathcal{L}_{\mathrm{HS}}(H, G)$ and that it satisfies

$$\mathbb{E}\bigg\{\int_0^T \|D_t X\|^2_{\mathcal{L}_{\mathrm{HS}}(H, G)}\, dt\bigg\} < +\infty$$

because of the growth condition (22) on the partial derivatives of $\kappa$ and the fact that Gaussian random variables have moments of all orders. It turns out that the Malliavin derivative $D$ as defined above as a densely defined operator from $L^2(\Omega; G)$ into $L^2([0, T] \times \Omega; \mathcal{L}_{\mathrm{HS}}(H, G))$ is closable. We use the same notation $D$ for its closure, and in particular, Definition 4.1 can be extended into the more practical one:

DEFINITION 4.2. If $X$ is the $L^2(\Omega, G)$ limit of a sequence $\{X_n\}_{n\geq 1}$ of smooth random variables such that $\{DX_n\}_{n\geq 1}$ converges in $L^2([0,T] \times \Omega; \mathcal{L}_{\mathrm{HS}}(H,G))$, we define
$$DX = \lim_{n \to \infty} DX_n.$$

REMARK 4.1 (Measurability). The Malliavin derivative $DX$ is defined to be an element of $L^2([0,T] \times \Omega; \mathcal{L}_{\mathrm{HS}}(H,G))$. Strictly speaking, it is an equivalence class of functions of $(t,\omega)$ which agree $\mathrm{Leb}_{[0,T]} \times \mathbb{Q}$ almost surely. By Fubini's theorem we can find a representative of $DX$ such that, for every $t \in [0,T]$, we have that $D_t X$ is measurable in $\omega$ and, for every $\omega \in \Omega$, we have that $DX(\omega)$ is measurable in $t$. We choose this representative to define $DX$.

We use the notation $\mathbb{H}^1(G)$ to represent the subspace of $L^2(\Omega; G)$ where the derivative can be defined by Definition 4.2. This subspace is a Hilbert space for the graph norm
$$\|X\|_{\mathbb{H}^1(G)}^2 = \mathbb{E}\{\|X\|_G^2\} + \mathbb{E}\left\{\int_0^T \|D_t X\|_{\mathcal{L}_{\mathrm{HS}}(H,G)}^2 \, dt\right\}.$$

The following simple sufficient condition for Malliavin differentiability will be needed in the sequel.

LEMMA 4.3. *If $X_n \to X$ converges in $L^2(\Omega; G)$, then we have $X \in \mathbb{H}^1(G)$ whenever the following boundedness condition is satisfied:*
$$\sup_n \mathbb{E}\left\{\int_0^T \|D_t X_n\|_{\mathcal{L}_{\mathrm{HS}}(H,G)}^2 \, dt\right\} < +\infty.$$

PROOF. The sequence $\{X_n\}$ is bounded in $\mathbb{H}^1(G)$, and hence, there exists a subsequence $\{X_{n_k}\}_k$ that converges weakly in $\mathbb{H}^1(G)$. But since $X_{n_k} \to X$ converges in $L^2(\Omega; G)$, we see that the weak limit of $\{X_{n_k}\}_k$ is $X$, implying that $X \in \mathbb{H}^1(G)$. □

Now we come to the Clark–Ocone formula, the crucial result that provides an explicit martingale representation for random variables in $\mathbb{H}^1(\mathbb{R})$ in terms of the Malliavin derivative. A version of this formula for stronger differentiability assumptions is originally due to Clark (1970). The formulation in terms of the Malliavin derivative is due to Ocone (1984).

THEOREM 4.4 (Clark–Ocone formula). *For every $\mathcal{F}_T$-measurable random variable $X \in \mathbb{H}^1(\mathbb{R})$, we have the representation*
$$X = \mathbb{E}\{X\} + \int_0^T \langle \mathbb{E}\{D_t X | \mathcal{F}_t\}, dW_t \rangle_H.$$

To prove this formula, we need the following integration by parts formula.

LEMMA 4.5. *Let $\{\beta_t\}_{t \in [0,T]}$ be an adapted process in $L^2([0,T] \times \Omega; H)$ and let $X \in \mathbb{H}^1(\mathbb{R})$. We have*

$$\mathbb{E}\left\{\int_0^T \langle D_t X, \beta_t\rangle_H \, dt\right\} = \mathbb{E}\left\{X \int_0^T \langle \beta_t^*, dW_t\rangle_H\right\}.$$

PROOF. First assume that $X = \kappa(\int_0^T \langle h_t^1, dW_t\rangle_H, \ldots, \int_0^T \langle h_t^n, dW_t\rangle_H)$ is a smooth random variable. Note that conditional on $\mathcal{F}_t$, the Wiener integral $\int_0^T \langle h_s, dW_s\rangle_H$ is a real Gaussian random variable with mean $\int_0^t \langle h_s, dW_s\rangle_H$ and variance $\int_t^T \|h_s\|_{H^*}^2 \, ds$, so that for $\mathcal{F}_t$-measurable $\beta \in L^2(\Omega; H)$, we have

$$\mathbb{E}\left\{\int_t^T \langle D_s X, \beta\rangle_H \, ds \,\Big|\, \mathcal{F}_t\right\} = \mathbb{E}\left\{X \int_t^T \langle \beta^*, dW_s\rangle_H \,\Big|\, \mathcal{F}_t\right\}$$

by definition of the Malliavin derivative for smooth random variables and the ordinary integration by parts formula.

Now assuming there exists $\mathcal{F}_{t_i}$-measurable $\beta_i \in L^2(\Omega; H)$ such that $\beta_t = \sum_{i=0}^N \mathbf{1}_{(t_i, t_{i+1}]}(t)\beta_i$, we have

$$\mathbb{E}\left\{\int_0^T \langle D_t X, \beta_t\rangle_H \, dt\right\} = \sum_{i=0}^N \mathbb{E}\left\{\mathbb{E}\left\{\int_{t_i}^{t_{i+1}} \langle D_t X, \beta_i\rangle_H \, dt \,\Big|\, \mathcal{F}_{t_i}\right\}\right\}$$

$$= \sum_{i=0}^N \mathbb{E}\left\{\mathbb{E}\left\{X \int_{t_i}^{t_{i+1}} \langle \beta_i^*, dW_t\rangle_H \, dt \,\Big|\, \mathcal{F}_{t_i}\right\}\right\}$$

$$= \mathbb{E}\left\{X \int_0^T \langle \beta_t^*, dW_t\rangle_H\right\}.$$

Since the smooth random variables are dense in $\mathbb{H}^1(\mathbb{R})$ and the simple integrands are dense in $L^2([0,T] \times \Omega; H^*)$, a straightforward limiting procedure completes the proof. $\square$

PROOF OF THEOREM 4.4. Since $X \in L^2(\Omega; \mathbb{R})$, by the martingale representation theorem there exists an adapted process $\{\alpha_t\}_{t \in [0,T]} \in L^2(\Omega \times [0,T]; H^*)$ such that

$$X = \mathbb{E}\{X\} + \int_0^T \langle \alpha_t, dW_t\rangle_H.$$

Without any loss of generality we can assume that $\mathbb{E}\{X\} = 0$. Now let $\{\beta_t\}_{t \in [0,T]}$ be an adapted measurable process in $L^2([0,T] \times \Omega; H)$. By Proposition 4.5 and Itô's isometry, we have

$$\mathbb{E}\left\{\int_0^T \langle D_t X, \beta_t\rangle_H \, dt\right\} = \mathbb{E}\left\{\int_0^T \langle \alpha_t, dW_t\rangle_H \int_0^T \langle \beta_t^*, dW_t\rangle_H\right\}$$

$$= \mathbb{E}\bigg\{\int_0^T \langle \alpha_t, \beta_t \rangle_H \, dt\bigg\},$$

implying

(23) $$\mathbb{E}\bigg\{\int_0^T \langle \lambda_t, \beta_t \rangle_H \, dt\bigg\} = 0,$$

where $\lambda_t = D_t X - \alpha_t$. The process $\{\lambda_t\}_{t \in [0,T]}$ is in $L^2([0,T] \times \Omega; H^*)$ by assumption, but it is not adapted to the filtration. Since the optional projection process $\{\mathbb{E}\{\lambda_t | \mathcal{F}_t\}\}_{t \in [0,T]}$ is obviously adapted to the filtration, and since $\{\lambda_t\}_{t \in [0,T]}$ is measurable and the filtration is right continuous, we have that $\{\mathbb{E}\{\lambda_t | \mathcal{F}_t\}\}_{t \in [0,T]}$ is adapted and measurable. Letting $\beta_t^* = \mathbb{E}\{\lambda_t | \mathcal{F}_t\}$ in (23), we get

$$\mathbb{E}\bigg\{\int_0^T \|\mathbb{E}\{\lambda_t | \mathcal{F}_t\}\|_{H^*}^2 \, dt\bigg\} = 0,$$

implying that $\alpha_t = \mathbb{E}\{D_t X | \mathcal{F}_t\}$ for almost every $(t, \omega)$ as desired. $\square$

We close this section with two results that allow us to calculate explicit formulas in what follows. The first one is a generalization of the chain rule in the spirit of Proposition 1.2.3 of Nualart (1995).

PROPOSITION 4.6. *Given a random variable $X \in \mathbb{H}^1(F)$ and a function $\kappa: F \to G$ such that*

$$\|\kappa(x) - \kappa(y)\|_G \leq C\|x - y\|_F$$

*for all $x, y \in F$ and some $C > 0$. Then $\kappa(X) \in \mathbb{H}^1(G)$ and there exists a random variable $\nabla \kappa(X)$ satisfying the bound $\|\nabla \kappa(X)\|_{\mathcal{L}(F,G)} \leq C$ almost surely and such that*

$$D\kappa(X) = \nabla \kappa(X) DX.$$

REMARK 4.2. We are not claiming that the function $\kappa$ is differentiable. Instead, we merely state that the random variable $\nabla \kappa(X)$ plays the role of a derivative in the sense of the chain rule. Of course if $\kappa$ is Fréchet differentiable, then $\nabla \kappa(X)$ is its Fréchet derivative evaluated at $X$. In Section 5, we use this result in the cases where $\kappa = g: F \to \mathbb{R}$ and when $\kappa = \sigma(t, \cdot): F \to \mathcal{L}_{\mathrm{HS}}(H, F)$.

PROOF OF PROPOSITION 4.6. According to Lemma 4.3, in order to show that $\kappa(X) \in \mathbb{H}^1(G)$, we need only find a sequence of functions $\{\kappa_n\}_n$ such that $\kappa_n(X) \to \kappa(X)$ strongly in $L^2(\Omega; F)$ and that $\{D\kappa_n(X)\}_n$ is bounded in $L^2([0,T] \times \Omega; H^*)$.

Let $\{e^i\}_{i=1}^{\infty}$ be a basis of $F$ and let $\{r^i\}_{i=1}^n$ be a basis for $\mathbb{R}^n$. Let

$$\ell_n = \sum_{i=1}^n e^i \otimes r^i \in \mathcal{L}(\mathbb{R}^n, F) \quad \text{and} \quad \ell'_n = \sum_{i=1}^n r^i \otimes e^{*i} \in \mathcal{L}(F, \mathbb{R}^n).$$

For every $n$, let $j_n : \mathbb{R}^n \to \mathbb{R}$ be a twice differentiable positive bounded function supported on the unit ball in $\mathbb{R}^n$ and such that $\int_{\mathbb{R}^n} j_n(x)\,dx = 1$, and for every $\varepsilon > 0$, define the approximate identity $j_n^\varepsilon$ by $j_n^\varepsilon(x) = \varepsilon^{-n} j_n(x/\varepsilon)$. Set $\varepsilon = 1/n$ and choose $\kappa_n$ to be defined by the Bochner integral

$$\kappa_n(x) = \int_{\mathbb{R}^n} j_n^\varepsilon(\ell'_n x - y)\kappa(\ell_n y)\,dy = \int_{\mathbb{R}^n} j_n^\varepsilon(y)\kappa(\ell_n \ell'_n x - \ell_n y)\,dy.$$

Note that $\kappa_n$ is differentiable and that

$$\mathbb{E}\{\|\kappa(X) - \kappa_n(X)\|_G^2\} \leq \mathbb{E}\bigg\{\bigg(\int_{\mathbb{R}^n} j_n^\varepsilon(y)\|\kappa(\ell_n \ell'_n X - \ell_n y) - \kappa(X)\|_G\,dy\bigg)^2\bigg\}$$

$$\leq C^2 \mathbb{E}\bigg\{\bigg(\int_{\mathbb{R}^n} j_n^\varepsilon(y)(\|(\ell_n \ell'_n X - X)\|_F + \|y\|_{\mathbb{R}^n})\,dy\bigg)^2\bigg\}$$

$$\leq 2C^2 \mathbb{E}\{\|(\ell_n \ell'_n - I)X\|_F^2\} + 2C^2/n^2 \to 0$$

by the dominated convergence theorem. By the finite-dimensional chain rule, we have

$$D\kappa_n(X) = \int_{\mathbb{R}^n} \kappa(\ell_n y) \otimes (\nabla j_n^\varepsilon(y - \ell'_n X) D\ell'_n X)\,dy,$$

where $\nabla$ is the gradient in $\mathbb{R}^n$, so that

$$\mathbb{E}\bigg\{\int_0^T \|D_t \kappa_n(X)\|_{\mathcal{L}_{\mathrm{HS}}(H,G)}^2\,dt\bigg\} \leq C^2 \mathbb{E}\bigg\{\int_0^T \|D_t X\|_{\mathcal{L}_{\mathrm{HS}}(H,F)}^2\,dt\bigg\},$$

and we can apply Lemma 4.3.

Finally, we note that $\nabla \kappa_n(X)$ is bounded in $L^\infty(\Omega; \mathcal{L}(F,G))$ and hence by the Banach–Alaoglu theorem there exists a subsequence $\{\nabla \kappa_{n_k}(X)\}_k$ and a random operator $\nabla \kappa(X)$ such that

$$\mathbb{E}\bigg\{\int_0^T \mathrm{trace}(A_t D_t \kappa_{n_k}(X))\,dt\bigg\} = \mathbb{E}\bigg\{\int_0^T \mathrm{trace}(A_t \nabla \kappa_{n_k}(X) D_t X)\,dt\bigg\}$$

$$\to \mathbb{E}\bigg\{\int_0^T \mathrm{trace}(A_t \nabla \kappa(X) D_t X)\,dt\bigg\}$$

for every $A \in L^2([0,T] \times \Omega; \mathcal{L}_{\mathrm{HS}}(F,H))$. On the other hand,

$$\mathbb{E}\bigg\{\int_0^T \mathrm{trace}(A_t D_t \kappa_{n_k}(X))\,dt\bigg\} \to \mathbb{E}\bigg\{\int_0^T \mathrm{trace}(A_t D_t \kappa(X))\,dt\bigg\}$$

so that $D_t \kappa(X) = \nabla \kappa(X) D_t X$ as claimed. □

The second result which we state without proof is the infinite-dimensional analog of (1.46) of Nualart (1995).

PROPOSITION 4.7. *If the adapted continuous square-integrable process* $\{\alpha_t\}_{t \in [0,T]}$ *is such that, for all* $t \in [0,T]$, *the random variable* $\alpha_t \in \mathbb{H}^1(\mathcal{L}_{\mathrm{HS}}(H,F))$ *is differentiable, then*

$$D_t \int_0^T \alpha_s \, dW_s = \alpha_t + \int_t^T D_t \alpha_s \, dW_s.$$

Note that when $\alpha$ and $W$ are scalar, the above result is true without assuming that $\alpha$ is adapted provided the stochastic integral is interpreted as a Skorohod integral instead of an Itô integral. We shall not need such a general form of this result.

**5. Hedging strategies for Lipschitz claims.** In this section we find explicit hedging strategies for an important class of contingent claims, and we characterize their properties. The results presented here are new. First we show that under natural conditions on the discounted bond price model and the payout function of the option, the hedging strategy is bounded in the $F^*$-norm, effectively avoiding the difficulties mentioned in Section 3 for hedging generic claims. Furthermore, we prove a general lemma which can be used to show that the hedging strategy is often confined to a small subspace of $F^*$. We apply this lemma to a model which has the essential features of a classical HJM model, yet exhibits some notion of maturity specific risk. For this class of models, we show that the counterintuitive strategies which are possible for finite-factor models are not allowed.

For the remainder of this article we make the following standing assumption.

ASSUMPTION 5.1. The contingent claim is European with expiration $T$ and payout given by $\xi = g(P_T)$. The payout function $g : F \to \mathbb{R}$ is such that the modified function $\tilde{g} : F \to \mathbb{R}$ given by $\tilde{g}(x) = x(T)g(x/x(T))$ satisfies the Lipschitz bound

$$|\tilde{g}(x) - \tilde{g}(y)| \le C_1 \|x - y\|_F \tag{24}$$

for all $x, y \in F$ and some constant $C_1 > 0$. Furthermore, for all $x_1, x_2 \in F$ such that $x_1(s) = x_2(s)$ for all $s \ge T$, we have

$$g(x_1) = g(x_2). \tag{25}$$

We remark that the condition (25) implies that the payout is insensitive to the part of the price curve corresponding to expired bonds. We also note that the Lipschitz assumption is reasonable. For instance, the payout function of a call option with expiration $T$ and strike $K$ on a bond with maturity $T_1 > T$ is $g(x) = (x(T_1) - K)^+$ and thus the modified payout function is given by

$$\tilde{g}(x) = x(T)\left(\frac{x(T_1)}{x(T)} - K\right)^+ = (x(T) - Kx(T))^+,$$

which is clearly Lipschitz since the point evaluations $\delta_T$ and $\delta_{T_1}$ are bounded linear functions on $F$.

If we can prove that $\tilde{P}_T \in \mathbb{H}^1(F)$, condition (24) and Proposition 4.6 imply that the Clark–Ocone formula applies. Our aim is to find an explicit representation of the Malliavin derivative $D\tilde{P}_T$ so that we can characterize the strategy that hedges $\xi$.

For the remainder of this article we work in a Markovian setting. The dynamics of the discounted bond prices will be given by Assumption 3.3 with the added provision that $\sigma_t = \sigma(t, \tilde{P}_t)$ for all $t \geq 0$. We list here the relevant assumptions on $\sigma(\cdot, \cdot)$.

ASSUMPTION 5.2. Let $\sigma(\cdot, \cdot) : \mathbb{R}_+ \times F \to \mathcal{L}_{\mathrm{HS}}(H, F_w^2)$ be such that $\sigma(\cdot, x)$ is continuous for all $x \in F$ such that $\sigma(t, 0) = 0$ for all $t > 0$, and such that we have the Lipschitz bound

$$\|\sigma(t, x) - \sigma(t, y)\|_{\mathcal{L}_{\mathrm{HS}}(H, F)} \leq C \|x - y\|_F \tag{26}$$

for all $t \geq 0$, $x, y \in F$ and some $C > 0$. We assume that, for every $x_1, x_2 \in F$ such that $x_1(s) = x_2(s)$ for all $s \geq t$, we have

$$\sigma(t, x_1) = \sigma(t, x_2). \tag{27}$$

REMARK 5.1. Recall we work in the space $F_w^2$ of differentiable functions described in Section 3 so that we can speak sensibly of interest rates, and in particular the bond price process is an Itô process. The conditions

$$\|\sigma(t, x)^* \delta_s\|_{H^*} \leq K |x(s)|,$$
$$\|\sigma(t, x) - \sigma(t, y)\|_{\mathcal{L}_{\mathrm{HS}}(H, F_w^2)} \leq K \|x - y\|_{F_w^2},$$

are sufficient to ensure that the bond prices are positive and condition (15) of Assumption 3.3 is satisfied. We will not make use of such conditions in the remainder of this paper. However, we are interested in hedging portfolios valued in the dual space $F^* = F_v^{1*}$, so we will explicitly make use of condition (26). Condition (27) implies that the volatility of the discounted prices is insensitive to the part of the curve corresponding to expired bonds.

First, we show that the Malliavin derivative of the discounted bond price exists.

LEMMA 5.3. *For all $T \geq 0$, we have that $\tilde{P}_T \in \mathbb{H}^1(F)$.*

PROOF. By Lemma 4.3, we need only to find a sequence of Malliavin differentiable random elements, say $\tilde{P}_T^n$, which converge toward $\tilde{P}_T$ in $L^2(\Omega; F)$, and such that $D\tilde{P}_T^n$ is bounded in $L^2([0,T] \times \Omega; \mathcal{L}_{\mathrm{HS}}(H,F))$. A natural candidate is provided by the elements of the Picard iteration scheme. Indeed, applying Proposition 4.7 to the $n$th step of the scheme, we obtain

$$D_t \tilde{P}_T^n = \sigma_t(\tilde{P}_t^{n-1}) + \int_t^T D_t \sigma_s(\tilde{P}_s^{n-1}) \, dW_s.$$

Now, since for all $s \in [0,T]$ the function $\sigma(s,\cdot)$ is Lipschitz, we can apply Proposition 4.6 and conclude by induction that $\tilde{P}_T \in \mathbb{H}^1(F)$. Indeed we have

$$\begin{aligned}
&\mathbb{E}\{\|D_t \tilde{P}_T^n\|_{\mathcal{L}_{\mathrm{HS}}(H,F)}^2\} \\
&= \mathbb{E}\{\|\sigma(t, \tilde{P}_t^{n-1})\|_{\mathcal{L}_{\mathrm{HS}}(H,F)}^2\} \\
&\quad + \mathbb{E}\left\{ \int_t^T \|D_t \sigma(s, \tilde{P}_s^{n-1})\|_{\mathcal{L}_{\mathrm{HS}}(H,\mathcal{L}_{\mathrm{HS}}(H,F))}^2 \, ds \right\} \\
&\leq C^2 \mathbb{E}\{\|\tilde{P}_t^{n-1}\|_F^2\} + C^2 \mathbb{E}\left\{ \int_t^T \|D_t \tilde{P}_s^{n-1}\|_{\mathcal{L}_{\mathrm{HS}}(H,F)}^2 \, ds \right\}.
\end{aligned}$$

Since the Picard iterates satisfy the bound

$$\mathbb{E}\{\|\tilde{P}_t^{n-1}\|_F^2\} \leq \|\tilde{P}_0\|_F^2 e^{C^2 t}$$

for all $n \geq 1$, we have

$$\mathbb{E}\{\|D_t \tilde{P}_T^n\|_{\mathcal{L}_{\mathrm{HS}}(H,F)}^2\} \leq C^2 \|\tilde{P}_0\|_F^2 e^{C^2 T},$$

by Gronwall's inequality. This completes the proof. $\square$

Since we know that $\tilde{P}_t \in \mathbb{H}^1(F)$ for all $t \geq 0$, we can conclude by Proposition 4.6 that for $\sigma(t, \tilde{P}_t) \in \mathbb{H}^1(\mathcal{L}_{\mathrm{HS}}(H,F))$, and by Proposition 4.7, we see that $\{D_t \tilde{P}_s\}_{s \in [t,T]}$ satisfies the linear equation

$$(28) \qquad D_t \tilde{P}_s = \sigma_t + \int_t^s \nabla \sigma_u D_t \tilde{P}_u \, dW_u.$$

Note that, for all $t \geq 0$, the random variable $\nabla \sigma_t$ takes values in $\mathcal{L}(F, \mathcal{L}_{\mathrm{HS}}(H,F))$, and that for each $s$, $t$, we have $D_t \tilde{P}_s \in \mathcal{L}_{\mathrm{HS}}(H,F)$, so that $\nabla \sigma_s D_t \tilde{P}_s \in \mathcal{L}_{\mathrm{HS}}(H, \mathcal{L}_{\mathrm{HS}}(H,F))$.

We now appeal to Skorohod's theory of strong random operators as developed in Skorohod (1984). A strong random operator from $F$ into $G$ is a $G$-valued stochastic process $\{Z_t(x)\}_{t \geq 0, x \in F}$ which is linear in $x \in F$. If such a process is adapted (in an obvious sense) and if, for example, $G = \mathcal{L}_{\mathrm{HS}}(H,F)$, then by setting

$$\left[ \int_a^b Z_s \, dW_s \right](x) = \int_a^b Z_s(x) \, dW_s,$$

we define a strong random operator $\int_a^b Z_s \, dW_s$ on $F$ (i.e., from $F$ into $F$). In particular, if for each $t > 0$, $\{Y_{t,s}(x)\}_{s \geq t, x \in F}$ is a strong random operator on $F$, then $\{\nabla \sigma_s Y_{t,s}(x)\}_{s \geq t, x \in F}$ is a strong random operator from $F$ into $\mathcal{L}_{\text{HS}}(H, F)$. Then with this definition of the integrand, the stochastic integral $\int_a^b \nabla \sigma_s Y_{t,s} \, dW_s$ is a strong operator on $F$. In this sense of equality of strong random operators on $F$, we would like to interpret the stochastic differential equation

$$(29) \qquad Y_{t,s} = I + \int_t^s \nabla \sigma_u Y_{t,u} \cdot dW_u,$$

where $I \in \mathcal{L}(F, F)$ is the identity. We are interested in solving such an equation because the solution process (if any) is in some sense the derivative of $\tilde{P}_T$ with respect to $\tilde{P}_t$. Moreover, the Malliavin derivative of the terminal underlying price should be related to this new process by

$$(30) \qquad D_t \tilde{P}_T = Y_{t,T} \sigma_t.$$

We settle the existence of a solution for this equation in the following proposition.

PROPOSITION 5.4. *Under the Lipschitz assumption* (26), *the linear equation*

$$Y_{t,s} = I + \int_t^s \nabla \sigma_u Y_{t,u} \cdot dW_u$$

*has a strong* $\mathcal{L}(F,F)$-*valued martingale solution* $\{Y_{t,s}\}_{s \in [t,T]}$. *Furthermore, we have the bound*

$$(31) \qquad \mathbb{E}\{\|Y_{t,T} x\|_F^2 | \mathcal{F}_t\} \leq \|x\|_F^2 e^{C^2(T-t)}.$$

PROOF. We prove that a Picard iteration scheme converges. Let $Y_{t,s}^0 = I$ for $s \in [t, T]$ and let

$$Y_{t,s}^{n+1} = I + \int_t^s \nabla \sigma_u Y_{t,u}^n \cdot dW_u.$$

Using the Lipschitz assumption and the martingale inequality, we have, for every $x \in F$,

$$\mathbb{E}\left\{\sup_{s \in [t,T]} \|(Y_{t,s}^{n+1} - Y_{t,s}^n)x\|_F^2 \Big| \mathcal{F}_t\right\}$$

$$\leq 4 \, \mathbb{E}\left\{\int_t^T \|\nabla \sigma_s (Y_{t,s}^n - Y_{t,s}^{n-1})x\|_{\mathcal{L}_{\text{HS}}(H,F)}^2 \, ds \Big| \mathcal{F}_t\right\}$$

$$\leq 4C^2 \int_t^T \mathbb{E}\{\|(Y_{t,s}^n - Y_{t,s}^{n-1})x\|_F^2 | \mathcal{F}_t\} \, ds.$$

So by induction, we have

$$\mathbb{E}\left\{\sup_{s\in[t,T]}\|(Y_{t,s}^{n+1}-Y_{t,s}^n)x\|_F^2\Big|\mathcal{F}_t\right\}\leq\|x\|_F^2\frac{C^{2n}(T-t)^n}{n!},$$

proving by a Borel–Cantelli lemma that the sequence of processes $\{Y_{t,\cdot}^n\}_n$ converges almost surely toward a process which is continuous in the strong topology of $\mathcal{L}(F,F)$. Furthermore, we have

$$\mathbb{E}\{\|Y_{t,s}x\|_F^2|\mathcal{F}_t\} = \|x\|_F^2 + \mathbb{E}\left\{\int_t^s \|\nabla\sigma_u Y_{t,u}x\|_{\mathcal{L}_{\mathrm{HS}}(H,F)}^2\, du\Big|\mathcal{F}_t\right\}$$

$$\leq \|x\|_F^2 + C^2\int_t^s \mathbb{E}\{\|Y_{t,u}x\|_F^2|\mathcal{F}_t\}\, du,$$

which implies the desired bound by Gronwall's inequality. $\square$

We assume that $\tilde{g}$ is Lipschitz, so we have the chain rule

$$D\tilde{g}(\tilde{P}_T) = \nabla\tilde{g}(\tilde{P}_T)D\tilde{P}_T,$$

where $D\tilde{P}_T \in \mathcal{L}_{\mathrm{HS}}(H,F)$ and $\nabla\tilde{g}(\tilde{P}_T) \in \mathcal{L}(F,\mathbb{R}) = F^*$.

We now use the Clark–Ocone formula, the chain rule and (30) to identify a candidate pre-hedging strategy from the following formal calculation:

$$\tilde{g}(\tilde{P}_T) = \mathbb{E}\{\tilde{g}(\tilde{P}_T)\} + \int_0^T \langle \mathbb{E}\{\nabla\tilde{g}(\tilde{P}_T)|\mathcal{F}_t\}, Y_{t,T}\sigma_t\, dW_t\rangle_F$$

$$= \mathbb{E}\{\tilde{g}(\tilde{P}_T)\} + \int_0^T \langle \mathbb{E}\{Y_{t,T}^*\nabla\tilde{g}(\tilde{P}_T)|\mathcal{F}_t\}, d\tilde{P}_t\rangle_F.$$

PROPOSITION 5.5.  *The process $\{\phi_t\}_{t\in[0,T]}$ given by the weak integral*

$$\phi_t = \mathbb{E}\{Y_{t,T}^*\nabla\tilde{g}(\tilde{P}_T)|\mathcal{F}_t\}$$

*is a pre-hedging strategy for the claim $g(P_T)$.*

REMARK 5.2.  By the formal calculation above, we only need to show that $\phi_t \in \mathcal{T}_{[t,\infty)}$ for almost all $(t,\omega) \in [0,T] \times \Omega$. First, we note that $\phi_t$ is bounded in $F^*$ uniformly in $t \in [0,T]$ almost surely. Indeed, we have

$$\|\phi_t\|_{F^*} = \sup_{\|x\|_F\leq 1} \mathbb{E}\{\langle\nabla\tilde{g}(\tilde{P}_T), Y_{t,T}x\rangle_F|\mathcal{F}_t\}$$

$$\leq \sup_{\|x\|_F\leq 1} C_1\mathbb{E}\{\|Y_{t,T}x\|_F^2|\mathcal{F}_t\}^{1/2}$$

$$\leq C_1 e^{C^2(T-t)/2}$$

by the Lipschitz bound (24) and the exponential growth bound (31). In fact, we have $\phi_t \in \mathcal{T}_{[t,\infty)}$ for almost all $(t,\omega) \in [0,T] \times \Omega$ thanks to the following lemma.

LEMMA 5.6. *Let $\{\mathcal{S}_t\}_{t\geq 0}$ be a decreasing family of closed subspaces of $F^*$ such that, for $s \leq t$, we have $\mathcal{S}_t \subset \mathcal{S}_s$ and such that $\delta_t \in \mathcal{S}_t$ for all $t \in [0,T]$. Suppose:*

(i) *the payout function $g(\cdot)\colon F \to \mathbb{R}$ is such that, for all $x \in F$ and all $y \in \mathcal{S}_T^\perp$, we have*

$$g(x+y) = g(x); \tag{32}$$

(ii) *the volatility function $\sigma(\cdot,\cdot)\colon \mathbb{R}_+ \times F \to \mathcal{L}_{\mathrm{HS}}(H,F)$ is such that, for each $t \geq 0$ and all $x \in F$, $y \in \mathcal{S}_t^\perp$ and $\mu \in \mathcal{S}_t$, we have*

$$\sigma(t, x+y)^*\mu = \sigma(t,x)^*\mu. \tag{33}$$

*Then for almost all $(t,\omega) \in [0,T] \times \Omega$, the random variable $\phi_t = \mathbb{E}\{Y_{t,T}^* \nabla \tilde{g}(\tilde{P}_T) | \mathcal{F}_t\}$ is valued in $\mathcal{S}_t$.*

PROOF. First we prove that

$$\langle \nabla \tilde{g}(X), y \rangle_F = 0 \quad \text{for any } X \in \mathbb{H}^1(F) \quad \text{and} \quad y \in \mathcal{S}_T^\perp, \tag{34}$$

where $\nabla \tilde{g}(X)$ is the bounded $F^*$-valued random variable such that $D\tilde{g}(X) = \nabla \tilde{g}(X) DX$. Following the proof of Proposition 4.6, we let $\tilde{g}_n(x)$ be given by

$$\tilde{g}_n(x) = \int_{\mathbb{R}^n} j_n^\varepsilon(u - \ell_n x) \tilde{g}(\ell_n u)\, du,$$

and recall that we have that $\tilde{g}_n(X)$ converges strongly to $\tilde{g}(X)$, and that there exists a subsequence such that $\nabla \tilde{g}_n(X)$ converges to $\nabla \tilde{g}(X)$ in the weak-* topology of $L^\infty(\Omega, F^*)$. Let $y \in \mathcal{S}_T^\perp$ and notice

$$|\langle \nabla \tilde{g}_n(X), y \rangle_F| = \left| \int_{\mathbb{R}^n} \langle \nabla j_n^\varepsilon(u - \ell_n x), \ell_n y \rangle_{\mathbb{R}^n} \tilde{g}(\ell_n u)\, du \right|$$

$$\leq \int_{\mathbb{R}^n} \lim_{h \to 0} \left| \frac{\tilde{g}(\ell_n' u + h \ell_n' \ell_n y) - \tilde{g}(\ell_n' u)}{h} \right| j_n^\varepsilon(u - \ell_n x)\, du$$

$$= \int_{\mathbb{R}^n} \lim_{h \to 0} \left| \frac{\tilde{g}(\ell_n' u + h \ell_n' \ell_n y) - \tilde{g}(\ell_n' u + hy)}{h} \right| j_n^\varepsilon(u - \ell_n x)\, du$$

$$\leq C_1 \|(I - \ell_n' \ell_n) y\|_F \to 0.$$

Similarly, we note that (33) implies that, for all $y \in \mathcal{S}_t^\perp$ and $\mu \in \mathcal{S}_t$, we have $(\nabla \sigma(t,X)y)^* \mu = 0$ and hence $\nabla \sigma(t,X) \in \mathcal{L}(\mathcal{S}_t^\perp, \mathcal{L}_{\mathrm{HS}}(H, \mathcal{S}_t^\perp))$.

The identity $I \in \mathcal{L}(F,F)$ obviously takes $\mathcal{S}_t^\perp$ into $\mathcal{S}_t^\perp$ and hence (29) has the strong operator solution $\{Y_{t,s}\}_{s \in [t,T]}$ valued in $\mathcal{L}(\mathcal{S}_t^\perp, \mathcal{S}_t^\perp)$. Thus for every $y \in \mathcal{S}_t^\perp$, we have

$$\langle \phi_t, y \rangle_F = \mathbb{E}\{\langle \nabla \tilde{g}(\tilde{P}_T), Y_{t,T} y \rangle_F | \mathcal{F}_t\} = 0$$

by (34), implying $\phi_t \in \mathcal{S}_t^{\perp\perp} = \mathcal{S}_t$ as desired. □

PROOF OF PROPOSITION 5.5. Apply Lemma 5.6 with the decreasing family of subspaces given by $\mathcal{S}_t = \mathcal{T}_{[t,\infty)}$. Note that the hypotheses are fulfilled since (27) implies

$$\sigma(t, x+y) = \sigma(t, x) \qquad \text{for all } y \in \mathcal{T}_{[t,\infty)}^{\perp},$$

and by assumption (25) we have $\tilde{g}(x+y) = \tilde{g}(x)$ for all $y \in \mathcal{T}_{[T,\infty)}^{\perp}$. □

Revisiting the motivating example of Section 3, for any contingent claim maturing at time $T$, we denote by $T' > T$ the longest maturity of the bonds underlying the claim. The following theorem shows that under the appropriate assumptions in the case of infinite-factor HJM models, the bonds in the hedging strategy for this claim have maturities less than or equal to $T'$. This intuitively appealing result is inspired by classical HJM models, of the type

$$d\tilde{P}_t(s) = \tilde{P}_t(s) \left\langle \int_t^s \kappa(f_t(u))\,du,\, dW_t \right\rangle_{\mathbb{R}^d}$$

for a deterministic function $\kappa : \mathbb{R} \to \mathbb{R}^d$. Note for these models, the volatility $\sigma_t^* \delta_s$ of the discounted bond price depends only on the forward rates $f_t(u) = -\frac{\tilde{P}_t'(u)}{\tilde{P}_t'(u)}$ for $u \in [t, s]$.

THEOREM 5.7. *Suppose that for every $s \geq t$, we have*

$$\sigma(t, x_1)^* \delta_s = \sigma(t, x_2)^* \delta_s$$

*whenever $x_1(u) = x_2(u)$ for all $u \in [t,s]$. If the payout function $g$ has the property that there exists a $T' > T$ such that $g(x_1) = g(x_2)$ for all $x_1, x_2 \in F$ such that $x_1(s) = x_2(s)$ for all $s \in [T, T']$, then there exists a hedging strategy $\{\varphi_t\}_{t \in [0,T]}$ that replicates the payout $g(P_T)$ and it is such that $\varphi_t \in \mathcal{T}_{[t,T']}$ for almost all $(t, \omega) \in [0, T] \times \Omega$.*

*Furthermore, if for all $x \in F$ and $t \geq 0$, we have*

$$\ker(\sigma(t,x)^*) = \mathcal{T}_{[0,t]},$$

*then the hedging strategy is unique.*

PROOF. Apply Lemma 5.6 with $\mathcal{S}_t = \mathcal{T}_{[t,T']}$ to the pre-hedging strategy given by

$$\phi_t = \mathbb{E}\{Y_{t,T}^* \nabla \tilde{g}(\tilde{P}_T) | \mathcal{F}_t\}.$$

Since $\phi_t \in \mathcal{T}_{[t,T']}$, we have the self-financing hedging strategy $\varphi_t = \phi_t + (V_t - \langle \phi_t, P_t \rangle_F) \delta_t$ is also valued in $\mathcal{T}_{[t,T']}$. Uniqueness follows from Proposition 3.7. □

This theorem implies that hedging strategies for this class of contingent claims have the property that the support of the portfolio at almost all times is confined to an interval. Moreover, the right endpoint of this interval is given by the longest maturity of the bonds underlying the claim, confirming our intuition about maturity specific risk.

**Acknowledgments.** The authors would like to thank the referee for a careful reading of the manuscript and for several suggestions which improved the presentation of the final version of the paper. Also we would like to mention that, after the completion of this work, we received a copy of the preprint of Ekeland and Taflin (2002) in which a related optimal bond portfolio problem is solved.

Department of Operations Research
 and Financial Engineering
Princeton University
Princeton, New Jersey 08544
USA
e-mail: rcarmona@princeton.edu

Department of Mathematics
University of Texas at Austin
Austin, Texas 78712
USA
e-mail: tehranch@math.utexas.edu